\documentclass[leqno]{article}

\usepackage{amsmath}
\usepackage{amssymb}
\usepackage{amsfonts}
\usepackage{enumerate}
\usepackage{vmargin}

\newcommand{\tens}{\otimes}

\newcommand{\N}{\ensuremath{\mathcal{N}}}
\newcommand{\B}{\ensuremath{\mathbb{B}}}

\newcommand{\M}{\ensuremath{\mathcal{M}}}

\newcommand{\R}{\ensuremath{\mathbb{R}}}
\newcommand{\C}{\ensuremath{\mathbb{C}}}

\newcommand{\la}{\langle}\newcommand{\ra}{\rangle}
\renewcommand{\le}{\ensuremath{\leqslant}}
\renewcommand{\ge}{\ensuremath{\geqslant}}
\renewcommand{\leq}{\ensuremath{\leqslant}}
\renewcommand{\geq}{\ensuremath{\geqslant}}
\newcommand{\n}{\noindent}
\newcommand{\8}{\infty}
\newcommand{\qed}{\hfill \vrule height6pt  width6pt depth0pt}

\renewcommand{\a}{\alpha}

\newcommand{\D}{\Delta}

\renewcommand{\th}{\theta}
\renewcommand{\l}{\lambda}

\newcommand{\f}{\varphi}

\newtheorem{thm}{Theorem}[section]

\newtheorem{prop}[thm]{Proposition}
\newtheorem{cor}[thm]{Corollary}
\newtheorem{lemma}[thm]{Lemma}

\newtheorem{remark}[thm]{Remark}

\newenvironment{rk}{\begin{remark}\rm}{\end{remark}}

\newenvironment{pf}[1][]{\noindent {\it Proof #1} : }{\hbox{~}\qed
\smallskip
}

\newcommand{\be}{\begin{eqnarray*}}
\newcommand{\ee}{\end{eqnarray*}}
\newcommand{\beq}{\begin{equation}}
\newcommand{\eeq}{\end{equation}}

\numberwithin{equation}{section}

\begin{document}

\title{Complex interpolation of weighted noncommutative $L_p$-spaces}
\date{}
\author{{\'E}ric Ricard and Quanhua Xu}

\maketitle

\begin{abstract}
 Let $\M$ be a semifinite von Neumann algebra equipped with a semifinite
normal faithful trace $\tau$. Let $d$ be an injective positive measurable operator with respect to $(\M,\,\tau)$ such that $d^{-1}$ is also measurable. Define
  $$L_p(d)=\left\{x\in L_0(\M)\;:\; dx+xd\in L_p(\M)\right\}\quad\mbox{and}\quad \|x\|_{L_p(d)}=\|dx+xd\|_p\,.$$
We show that for $1\le p_0<p_1\le\8$, $0<\th<1$ and $\a_0\ge0, \a_1\ge0$ the interpolation equality
 $$(L_{p_0}(d^{\a_0}),\;L_{p_1}(d^{\a_1}))_\th
 =L_{p}(d^{\a})$$
holds with equivalent norms, where $\frac1p=\frac{1-\th}{p_0}+\frac{\th}{p_1}$ and $\a=(1-\th)\a_0+\th\a_1$.
 \end{abstract}


\makeatletter
 \renewcommand{\@makefntext}[1]{#1}
 \makeatother \footnotetext{\noindent
 This work is partially supported by ANR 06-BLAN-0015.
 }


\section{Introduction}


Let $\M$ be a semifinite von Neumann algebra equipped with a
normal faithful semifinite trace $\tau$. For $1\le p\le\8$, let
$L_p(\M)$ denote the noncommutative $L_p$-space associated with
$(\M,\,\tau)$. The norm of $L_p(\M)$ is denoted by $\|\;\|_p$. All
spaces $L_p(\M)$ are continuously injected into the topological
involutive algebra $L_0(\M)$ of measurable operators with respect
to $(\M, \,\tau)$. This injection turns
$(L_{p_0}(\M),\;L_{p_1}(\M))$ into a compatible couple. We
then have the following well-known identity on the complex
interpolation of noncommutative $L_p$-spaces: for any $1\le
p_0, p_1\le\8$ and $0<\th<1$,
 \beq\label{inter}
 (L_{p_0}(\M),\;L_{p_1}(\M))_\th=L_{p}(\M)
 \eeq
with equal norms, where
$\frac1p=\frac{1-\th}{p_0}+\frac{\th}{p_1}$. We refer to
\cite{fk}, \cite{nelson}, \cite{terp} and \cite{px-survey} for
semifinite noncommutative $L_p$-spaces and to \cite{bl} for
interpolation.

The aim of this note is to consider the weighted version of
\eqref{inter}. Let $d\in L_0(\M)$ be a positive injective operator
such that $d^{-1}\in L_0(\M)$. We will call $d$ a density.
Define
 $$L_p^r(d)=\left\{x\in L_0(\M)\;:\; xd\in L_p(\M)\right\}$$
equipped with the norm
 $$\|x\|_{L_p^r(d)}=\|xd\|_p\,.$$
Then by standard arguments one easily deduces from \eqref{inter}
the following right weighted analogue. Let $\th, p_0, p_1, p$
be as in \eqref{inter}, and let $\a_0, \a_1\in\R$,
$\a=(1-\th)\a_0+\th\a_1$. Then
 \beq\label{winter1}
 (L_{p_0}^r(d^{\a_0}),\;L_{p_1}^r(d^{\a_1}))_\th
 =L_{p}^r(d^{\a})
 \eeq with equal norms. One has, of course, a similar equality for
the left weighted spaces. However, the matter becomes highly
subtle when one considers the sum of two multiplication maps, one
from left and another from right. Thus let
 $$L_p(d)=\left\{x\in L_0(\M)\;:\; dx+xd\in L_p(\M)\right\}$$
equipped with the norm
 $$\|x\|_{L_p(d)}=\|dx+xd\|_p\,.$$
We will see later that $L_p(d)$ is complete, so is a Banach space
for any $1\le p\le\8$. The compatibility on these weighted spaces
is induced by the identity of $L_0(\M)$. The following is the
two-sided version of \eqref{winter1}, which is the main result of
this note.

\begin{thm}\label{winter11}
 Let $0<\th<1$, $1\le p_0,p_1\le\8$ and
$\frac1p=\frac{1-\th}{p_0}+\frac{\th}{p_1}$.
Assume $1<p<\8$. Let $\a_0, \a_1\ge0$
and $\a=(1-\th)\a_0+\th\a_1$.
Then
 \beq\label{wi1-1}
 (L_{p_0}(d^{\a_0}),\;L_{p_1}(d^{\a_1}))_\th
 =L_{p}(d^{\a})
 \eeq
with equivalent norms.
 \end{thm}

It is worth to note that while \eqref{winter1} holds for any real $\a_0$ and $\a_1$,
\eqref{wi1-1} may fail when $\a_0$ and $\a_1$ are of opposite signs (see Remark~\ref{os} below).

The theorem above is closely related to a recent remarkable
interpolation theorem of Junge and Parcet \cite{jupar-ros}. Define
  $$\D_p(d)=\left\{x\in L_0(\M)\;:\; dx,\,xd\in L_p(\M)\right\}$$
and
 $$\|x\|_{\D_p(d)}=\max\left(\|dx\|_p\,,\;\|xd\|_p\right)\,.$$
Keeping the assumptions of Theorem~\ref{winter11} on the indices,
we have
 \be
 (\D_{p_0}(d^{\a_0}),\;\D_{p_1}(d^{\a_1}))_\th
 =\D_{p}(d^{\a})
 \ee with equivalent norms. This is
\cite[Theorem~1.15]{jupar-ros}. Note that the von Neumann algebra
there should be assumed semifinite. Theorem~\ref{winter11} can be
obtained by duality from Junge and Parcet's theorem. We prefer,
however, to give a direct proof. In fact, we will show a slightly
stronger result (see Theorem~\ref{wi11} below). Note that the
pattern of our arguments still models that of \cite{jupar-ros}.
The two main ingredients are again the boundedness of some special
Schur multipliers and Pisier's interpolation theorem on triangular
subspaces of Schatten classes. Thus our arguments are very similar
to those of \cite{jupar-ros}.

Using Haagerup's reduction theorem (see \cite{hjx}), we deduce
from Theorem~\ref{winter11} a similar result for type III von
Neumann algebras as in \cite{jupar-ros}. Let $\N$ be a von Neumann
algebra equipped with a normal faithful state $\psi$. Let
$L_p(\N)$ be the Haagerup noncommutative $L_p$-spaces associated
with $\N$ (see \cite{haag-Lp} and \cite{terp}). Recall that
$L_\8(\N)=\N$ and $L_1(\N)$ coincides with the predual of $\N$.
Let $d$ be the operator in $L_1(\N)$ corresponding to $\psi$. For
$1\le p<\8$, consider the injection
 $$\iota_p\,:\, \N\to L_p(\N),\quad
 \iota_p(x)=d^{1/p}x+xd^{1/p}\,.$$ Note that $\iota_p$ is
injective and of dense range. $\iota_1$ makes $(\N,\; L_1(\N))$
into a compatible couple. We then have the following two-sided
analogue of Kosaki's  interpolation theorem \cite{kos-int}.

\begin{cor}\label{wih}
 Let $1<p<\8$.
Then
 $$
 (\N,\; L_1(\N))_{1/p}
 =L_{p}(\N)
 $$
with equivalent norms. More precisely, there exist two positive
constants $c_p$ and $C_p$ depending only on $p$ such that
 $$c_p\|d^{1/p}x+xd^{1/p}\|_p\le
 \|x\|_{(\N,\; L_1(\N))_{1/p}}
 \le C_p\|d^{1/p}x+xd^{1/p}\|_p\,,\quad \forall\; x\in\N\,.$$
 \end{cor}


\section{Schur multipliers}


In this section we consider some special Schur multipliers on
$\B(\ell_2)$, which will play a key role in the proof of our
interpolation theorem. These multipliers are of the type already discussed in \cite{jupar-ros}. Our presentation is, however, independent of \cite{jupar-ros}. As usual, the operators in $\B(\ell_2)$ are
represented as infinite matrices $x=(x_{ij})_{i,j\ge1}$ (with
respect to the canonical matrix units $\{e_{ij}\}$). Recall that a
bounded Schur multiplier on $\B(\ell_2)$ is an infinite matrix
$\f=(\f_{ij})_{i,j\ge1}$ of complex numbers such that
$(\f_{ij}x_{ij})_{i,j\ge1}\in \B(\ell_2)$ for any
$(x_{ij})_{i,j\ge1}\in \B(\ell_2)$. The resulting bounded map
$(x_{ij})_{i,j}\mapsto(\f_{ij}x_{ij})_{i,j}$ will be also denoted
by $\f$.

It is well known that $\f$ is a bounded Schur multiplier on
$\B(\ell_2)$ iff there exists a Hilbert space $H$ and two bounded
sequences $(\xi_i), (\eta_i)\subset H$ such that
 $$\f_{ij}=\la \xi_i\,,\; \eta_j\ra,\quad\forall\; i, j\ge1$$ (see
\cite[Theorem~5.1]{pis-sim}). Moreover, in this case $\f$ is
automatically completely bounded and $\|\f\|=\|\f\|_{cb}$. Recall that
the cb-norm $\|\f\|_{cb}$ of $\f$ is the norm of the map
$\rm{id}_{\B(\ell_2)}\tens\f$ on $\B(\ell_2)\bar\tens\B(\ell_2)$
(usually one uses the compact operators $\mathbb K(\ell_2)$ instead of
$\B(\ell_2)$ but it does not make any difference).  Let
$M_{cb}(\B(\ell_2))$, or simply $M_{cb}$ denote the space of all
completely bounded Schur multipliers $\f$ on $\B(\ell_2)$, equipped
with the norm $\|\f\|_{M_{cb}}$.  We then have
 $$\|\f\|_{M_{cb}}=\inf\big\{\sup_{i,j}\|\xi_i\|\,\|\eta_j\|\;:\;
 \f_{ij}=\la \xi_i\,,\; \eta_j\ra,\; \xi_i, \eta_j\in H,\; H\;
 \textrm{a Hilbert space}\big\}.$$

The following is a well-known elementary fact (see
\cite{jupar-ros} for a similar statement with more regularity on
the function $f$). We include a proof for the convenience of the
reader. As usual, $\hat f$ denotes the Fourier transform of a
function $f\in L_1(\R)$:
 $$\hat f(\xi)=\int_{\R} f(s) e^{-2\pi{\rm i}\xi s} {\rm d}s.$$

\begin{prop}\label{mult1}
 Let $f\in L_1(\R)$ such that  $\hat f$
belongs to $L_1(\R)$. Then, for any $s_i\in \R$,
$(f(s_i-s_j))_{i,j}\in M_{cb}$ and
 $$\|(f(s_i-s_j))_{i,j}\|_{M_{cb}}
 \leq \|\hat f\|_1.$$
In particular, if $f: \R_+\to \R_+$ is a nonincreasing convex
function, then, for any $s_i\in \R$, $(f(|s_i-s_j|))_{i,j}$ defines
a completely positive Schur multiplier on $\B(\ell_2)$.
\end{prop}

\begin{pf}
By the Fourier inversion formula, we have
 $$f(s)=  \int_{\R} \hat f(\xi) e^{{2\pi}{\rm i}\xi s} {\rm d}\xi\,.$$
Thus letting $g_i(\xi)=\sqrt{|\hat f(\xi)|}\, e^{{\rm i}s_i\xi}$
and $h_j=\overline {\hat f(\xi)} /\sqrt{|\hat f(\xi)|}\, e^{{\rm
i}s_j\xi}$, we have
 $$\|g_i\|_2^2=\|h_j\|_2^2=\|\hat f\|_1\quad\textrm{and}\quad
 f(s_i-s_j)=\la g_i,\;h_j\ra_{ L_2(\R)}\,;$$
whence the first assertion.

For the second one, we can suppose $\lim_{t\to+\infty} f(t)=0$
since the constant matrix $(1)_{i,j}$ is a completely positive
Schur multiplier on $\B(\ell_2)$. Put $g(s)=f(|s|)$ and note that
we do not need to assume that $g\in L_1(\R)$ to define the Fourier
transform (except for $\xi=0$).  Then, it is well known that $g$
is a positive definite function on $\R$, that is, $\hat g\geq 0$
and $\|\hat g\|_1=g(0)=f(0)$, and moreover the Fourier inversion
formula holds for $g$, see \cite[Theorem~8.7]{cha} for instance.
 \end{pf}

\begin{cor}\label{correc0}
 Let $(\lambda_i)_{i\geq 1}$ be a sequence of positive real numbers. Then,
for any $\theta\in [0,\,1]$,
 $$\begin{array}{ll}
 \displaystyle
 \left\|\left(\frac{\min(\lambda_i,\,\lambda_j)}
 {\max(\lambda_i,\,\lambda_j)}\right)_{i,j}\right\|_{M_{cb}}\le1\,,\quad
 &\displaystyle\left\|\left(\frac{(\lambda_i+\lambda_j)^\th}
 {\max(\lambda_i,\,\lambda_j)^\th}\right)_{i,j}\right\|_{M_{cb}}\leq2^\th\,,\\
 \displaystyle \left\|\left(\frac
 {\max(\lambda_i,\,\lambda_j)^\theta}
 {(\lambda_i+\lambda_j)^\theta}\right)_{i,j}\right\|_{M_{cb}}
 \leq 2- 2^{-\theta}\,,
 &\displaystyle \left\|\left(\frac{\min(\lambda_i,\,\lambda_j)^\theta}
 {(\lambda_i+\lambda_j)^\theta}\right)_{i,j}\right\|_{M_{cb}}
 \leq2-2^{-\theta}\,.
 \end{array}$$
 \end{cor}

\begin{pf}
 The first multiplier is obtained by applying Proposition~\ref{mult1} to the decreasing convex function $f(t)=e^{-t}$. Similarly,  the function corresponding to the second multiplier is
$f(t)=(1+e^{-t})^\th$. To deal with the third one, notice that with $\lambda_i=e^{s_i}$, we have
 $$\frac{\max(\lambda_i,\,\lambda_j)^\theta}{(\lambda_i+\lambda_j)^\theta}
 =\frac 1 {(1+e^{-|s_i-s_j|})^{\theta}}=-\Big(1-\frac 1
 {(1+e^{-|s_i-s_j|})^{\theta}}\Big)+1.$$
The function $f(t)=1-\frac 1{(1+e^{-t})^{\theta}}$ is decreasing
and convex on $\R_+$, so we get the estimate for the third  multiplier. The last one
is just the composition of the third with
$\big(\frac{\min(\lambda_i^\theta,\,\lambda_j^\theta)}
{\max(\lambda_i^\theta,\,\lambda_j ^\theta)}\big)_{i,j}$, which is
a complete contraction.
\end{pf}

\begin{cor}\label{correc}
 Let $(\lambda_i)_{i\geq 1}$ and $(\mu_i)_{i\geq 1}$
be two nondecreasing  sequences of positive real numbers. Then, for any
$\theta\in [0,\,1]$,
 $$\left\|\left(\frac {\lambda_i^{1-\theta}\, \mu_i^\theta+
 \lambda_j^{1-\theta}\, \mu_j^\theta}
 {(\lambda_i+\lambda_j)^{1-\theta}(\mu_i+\mu_j)^{\theta}}
 \right)_{i,j}\right\|_{M_{cb}}\leq 9-4\sqrt 2\,.$$
 \end{cor}

\begin{pf}
 Using
 $$\lambda_i^{1-\theta} \mu_i^\theta+
 \lambda_j^{1-\theta} \mu_j^\theta=
 \max(\lambda_i,\,\lambda_j)^{1-\theta}\max(\mu_i,\,\mu_j)^{\theta}
 +\min(\lambda_i,\,\lambda_j)^{1-\theta}\min(\mu_i,\,\mu_j)^{\theta}\,,$$
we immediately get the estimate from Corollary \ref{correc0}.
 \end{pf}

\begin{cor}\label{correc2}
 Let $(\lambda_i)_{i\geq 1}$ and $(\mu_i)_{i\geq 1}$ be nondecreasing
sequences of positive real numbers. Then, for any $\theta\in [0,\,1]$,
 $$\left\|\left(\frac{\big(\lambda_i+
 \lambda_j\big)^\theta\big(\mu_i+\mu_j\big)^{1-\theta}}
 {\lambda_i^\theta\mu_i^{1-\theta}+\lambda_j^\theta\mu_j^{1-\theta}}
 \right)_{i,j}\right\|_{M_{cb}}\leq 3.$$
\end{cor}

\begin{pf}
 Since $(\lambda_i)_{i\geq 1}$ and $(\mu_i)_{i\geq 1}$
have the same variation, we can write the multiplier under
consideration as a composition of three multipliers:
 $$\frac {\big(\lambda_i+\lambda_j\big)^\theta}
 {\max(\lambda_i^\theta,\,\lambda_j^\theta)}\,
 \frac{\big(\mu_i+\mu_j\big)^{1-\theta}}
 {\max(\mu_i^{1-\theta},\,\mu_j^{1-\theta})}\,
 \frac{\max(\lambda_i^\theta\mu_i^{1-\theta},\,
 \lambda_j^\theta\mu_j^{1-\theta})}
 {\lambda_i^\theta\mu_i^{1-\theta}+\lambda_j^\theta\mu_j^{1-\theta}}\,.$$
Then the assertion follows from Corollary \ref{correc0}.
\end{pf}

\begin{cor}\label{g-mean}
 Let $(\lambda_i)_{i\geq 1}$ be a sequence of positive real numbers. Then,
for any $0<\theta<1$,
 $$\left\|\left(\frac
 {\lambda_i^\theta\,\lambda_j^{1-\theta}}
 {\lambda_i+\lambda_j}\right)_{i,j}\right\|_{M_{cb}}\leq C
 \ln \frac 1{\theta(1-\theta)}\,.$$
Here, as well as in the sequel, $C$ denotes a universal positive constant.
 \end{cor}

\begin{pf}
According to Proposition~\ref{mult1}, we have to compute the
$L_1$-norm of the Fourier transform  of $f(s)=(e^{\theta
s}+e^{(\theta-1)s})^{-1}$. A standard calculation by the residue
theorem yields
 $$\hat f (\xi)=
 \frac \pi {\sin \big( \pi( \theta+ 2{\rm i}\pi\xi)\big)}\,.$$
Then it remains to note that $\|\hat f\|_1$ behaves like $\int_0^1
|\theta + {\rm i} x|^{-1} {\rm d} x$ when $\theta$ is close to
$0$.
\end{pf}

\begin{rk}
 The preceding corollary is to be compared with
 \cite[Lemma~1.7]{jupar-ros}, which asserts that $\big(\frac
 {\lambda_i^\theta\,\lambda_j^{1-\theta}}
 {\lambda_i+\lambda_j}\big)_{i,j}$ is a bounded  multiplier on the
 triangular subalgebra of $\B(\ell_2)$, for any $0\le\th\le 1$.
 \end{rk}

\begin{rk}\label{trans}
 Let $S_p$ denote the Schatten $p$-class,
i.e., $S_p=L_p(\B(\ell_2))$ with $\B(\ell_2)$ equipped with the
usual trace.  Similarly, we define  Schur multipliers on $S_p$ as
before for $\B(\ell_2)$. It is well-known that any bounded Schur
multiplier on $\B(\ell_2)$ is also bounded (even completely
bounded) on $S_p$ for any $1\le p<\8$.  Let us state this in a
slightly more general setting that will be crucial for the next
section. Let $\f$ be a bounded Schur multiplier on $\B(\ell_2)$
and $\M$ a von Neumann algebra. Then ${\rm id}_{L_p(\M)}\tens\f$
defines a bounded map on $L_p(\M\bar\tens\B(\ell_2))$, for any
$1\le p\le\8$:
 $$\big\|\left(\f_{ij}x_{ij}\right)_{i,j}
 \big\|_{L_p(\M\bar\tens\B(\ell_2))}
 \le\|\f\|_{cb}\,\big\|\left(x_{ij}\right)_{i,j}
 \big\|_{L_p(\M\bar\tens\B(\ell_2))}$$
for all finite matrices $(x_{ij})$ with entries in $L_p(\M)$.
 \end{rk}


\section{Interpolation}


In this section $(\M,\tau)$ will denote a semifinite von Neumann
algebra and  $d$ a density in $L_0(\M)$ such that $d^{-1}\in
L_0(\M)$. For $1\le p\le\8$, we define
 $$L_p(d)=\left\{x\in L_0(\M)\;:\; dx+xd\in L_p(\M)\right\}
 \quad\textrm{and}\quad \|x\|_{L_p(d)}=\|dx+xd\|_p\,.$$
Then $L_p(d)$ is a Banach space. The nontrivial point is the
completeness of the norm. This is an immediate consequence of
Proposition~\ref{inv} below.

We will use $L_d$ to denote the left multiplication map by $d$,
i.e., $L_d(x)=dx$. Similarly, $R_d$ is the right multiplication
map by $d$. It is clear that both $L_d$ and $R_d$ are continuous
on $L_0(\M)$. We will also consider them as closed densely defined
maps on $L_2(\M)$. In this latter case, they are injective,
positive and commuting. Thus we can apply functional calculus to
them. In particular, $\frac{\sqrt{L_dR_d}}{L_d+R_d}$ is also an
injective positive map on $L_2(\M)$.

 \begin{prop}\label{inv}
 $\frac{\sqrt{L_dR_d}}{L_d+R_d}$ extends to a bounded map
on $L_p(\M)$, for any $1\le p\le\8$, with norm
$\,\le1/2$. More precisely, we have the following integral
representation
 \beq\label{sqr}
 \frac{\sqrt{L_dR_d}}{L_d+R_d}\,(x)=
 \int_\R u_t(x)\,\frac{{\rm d}t}{2\cosh(\pi t)}\,,
 \eeq
where $(u_t)_{t\in\R}$ is the isometry group on $L_p(\M)$ defined
by
 $$u_t(x)=e^{{\rm i}t\ln d}xe^{-{\rm i}t\ln d}\,.$$
Consequently, $(L_d+R_d)^{-1}$ is a continuous map from $L_p(\M)$
to $L_0(\M)$ and $L_d+R_d$ is a isometry from $L_p(d)$ onto $L_p(\M)$.
 \end{prop}

\begin{pf} Consider first the case where
$d=\sum_{i=1}^k \lambda_i e_i$ for some increasing sequence
$(\lambda_i)$ of positive real numbers and mutually orthogonal
projections $e_i$ with sum $1$. Note that in this case $L_d$ and
$R_d$ are bijections on $L_p(\M)$, for any $1\le p\le\8$. It is
also clear that
 $$L_d(x)=\sum_i\lambda_i e_ix \quad\textrm{and}\quad
 R_d(x)=\sum_j\lambda_j xe_j\,.$$
Then $\frac{\sqrt{L_dR_d}}{L_d+R_d}$ is given by
 $$\frac{\sqrt{L_dR_d}}{L_d+R_d}\,(x)=\sum_{i,j=1}^k
 \frac{\sqrt{\lambda_i \lambda_j}}{\lambda_i+\lambda_j}\,
 e_ixe_j\,,\quad \forall\; x\in L_p(\M).$$
Applying Corollary~\ref{g-mean} (or its proof) with $\th=1/2$ we
find
 $$\frac{\sqrt{\lambda_i \lambda_j}}{\lambda_i+\lambda_j}
 =\int_\R e^{{\rm i}t(\ln \l_i-\ln \l_j)}\,\frac{{\rm d}t}{2\cosh(\pi t)}\,.$$
Thus \eqref{sqr} follows and
 $$\left\|\frac{\sqrt{L_dR_d}}{L_d+R_d}\right\|_{\B(L_p(\M))}
 \le\frac12\,.$$
We then deduce the general case by a standard approximation argument
(see also  step~2 of the proof of Theorem~\ref{wi11} below).

For the second part we note that
 $$(L_d+R_d)^{-1}=L_d^{-1/2}\,\frac{\sqrt{L_dR_d}}{L_d+R_d}
 \,R_d^{-1/2}\,.$$
Since the multiplication map by $d^{-1/2}$ from both left and
right is continuous from $L_p(\M)$ to $L_0(\M)$, we obtain the
desired continuity of $(L_d+R_d)^{-1}$.
\end{pf}

\begin{thm}\label{wi11}
 Let $f_0, f_1:\R_+\to\R_+$ be two  nondecreasing functions
with $f_i(t)>0$, for $t>0$. Put $d_0=f_0(d)$ and $d_1=f_1(d)$. Let
$1\leq p_0, p_1\leq \infty$ and $0<\theta<1$. Set
 $$\frac 1 {p}=\frac {1-\theta} {p_0}+\frac {\theta} {p_1}
 \quad\mbox{and}\quad
 d_\theta=d_0^{1-\theta}d_1^\theta\,.$$
Assume $1<p<\8$. Consider $(L_{p_0}(d_0),L_{p_1}(d_1))$ as a
compatible couple by injecting both spaces into $L_0(\M)$. Then,
for $0<\theta<1$, we have
 $$\left(L_{p_0}(d_0),\;L_{p_1}(d_1)\right)_\theta =L_{p}(d_\theta)\,.$$
More precisely, for $x\in L_{p_0}(d_0)\cap L_{p_1}(d_1)$
 \beq\label{int-ineq}
 C_{p'}^{-1}\|x\|_{L_{p}(d_{\th})}\le
 \|x\|_{(L_{p_0}(d_0),\;L_{p_1}(d_1))_\theta}
 \le C_p\|x\|_{L_{p}(d_{\th})}\,,
 \eeq
where $p'$ denotes the conjugate index of $p$, and where the
constant $C_p$ satisfies the following estimate
 $$C_p\le C\max(p,\, 2)\max(p,\, p')\,.$$
\end{thm}

The proof of Theorem~\ref{wi11} will be divided into two steps.
The first one deals with the case where $d$ has only point
spectrum. This is the main step. The second one is a simple
approximation argument.

\medskip\n{\it Step~1: The discrete case.} We assume that
$d=\sum_{i=1}^k \lambda_i e_i$ for some increasing sequence
$(\lambda_i)$ of positive real numbers and mutually orthogonal
projections $e_i$ with sum $1$. Then the map
 $$\kappa: \sum_{i,j} e_i xe_j \mapsto \sum_{i,j} e_i xe_j\tens e_{ij}$$
defines an isometry from $L_p(\M)$ into
$L_p(\M\bar\tens\B(\ell_2))$, for any $1\le p\le\8$. Moreover, its
range is contractively complemented. We will need the triangular
projections:
 $$T_+(x)=\sum_{j\geq i} e_ixe_j
 \quad\mbox{and}\quad T_-(x)=x-T_+(x).$$
Both $T_+$ and $T_-$  commute with $L_{f(d)}$ and $R_{f(d)}$.

\begin{lemma}\label{compa}
 For any $1\leq p\leq \infty$ and $f:\R_+\to \R_+$ nondecreasing,
 $$\frac23\, \|T_\pm(x)f(d)\|_{p}\leq \|T_\pm(x)\|_{L_p(f(d))}\leq 2
 \|T_\pm(x)f(d)\|_{p}\,.$$
\end{lemma}

\begin{pf}
  We have
 \begin{eqnarray*}
 &T_+(x)f(d) =\sum_{j\geq i} f(\lambda_j) e_ixe_j=
 \sum_{j\geq i} \max(f(\lambda_i),\,f(\lambda_j)) e_ixe_j\,, \\
 & f(d) T_+(x)+T_+(x)f(d)= \sum_{j\geq i}
 (f(\lambda_i)+f(\lambda_j)) e_ixe_j\,.&
 \end{eqnarray*}
As for any Schur multiplier $\f=(\f_{ij})$, we have the identity
 $$\kappa\big(\sum_{i,j} \f_{ij} e_i xe_j\big)=
 ({\rm id}\tens\f)\kappa\big(\sum_{i,j} e_i xe_j\big),$$
we then deduce the estimates on $T_+$ using Corollary~\ref{correc0}  and the transference
principle in Remark~\ref{trans}.
\end{pf}

\begin{rk}\label{triproj}
The transference principle also shows that  the triangular
projections are bounded on $L_p(\M)$, for any $1<p<\infty$:
 $$\|T_\pm(x)\|_p \leq C  \max(p,\,p') \|x\|_p\,,\quad
 \forall\;x\in L_p(\M)$$
for the norms of the triangular projections on $S_p$ are of order
$\max(p,\,p')$.
\end{rk}

We will use Pisier's interpolation theorem  on subspaces of
triangular matrices in \cite{pis-intHpI}. Let $T_p(\M)\subset
L_p(\M\bar\tens\B(\ell_2))$ be the subspace of upper triangular
matrices with respect to the matrix units $\{e_{ij}\}$ of
$\B(\ell_2)$. For any $p_0$ and $p_1$, the couple $(T_{p_0}(\M),\;
T_{p_1}(\M))$ is compatible in a natural way.

\begin{lemma}\label{pisierinter}\rm{({\bf Pisier})}
 We have
 $$(T_{p_0}(\M),\; T_{p_1}(\M))_\theta =T_{p}(\M)$$
with equivalent norms. More precisely, for any $x\in
T_{p_0}(\M)\cap T_{p_1}(\M)$
 $$\|x\|_p\le \|x\|_{(T_{p_0}(\M),\; T_{p_1}(\M))_\theta}
 \le t_p\|x\|_p\,,$$
where the constant $t_p$ is estimated by
 $t_p\le C\max(p,\; 2).$
\end{lemma}

\begin{rk}
 The estimate above of $t_p$ is not explicitly stated in
\cite{pis-intHpI}. It can be, however, tracked from Pisier's
proof. First, by \cite[Theorem~3.3.1]{bl} and
\cite[Theorem~4.3]{holm}, we find
 $$t_p\le C\max(p,\; p').$$
Next, to see that $t_p$ remains bounded when $p\to1$, we use the
reiteration theorem and the fact that
 $$(T_{1}(\M),\; T_{2}(\M))_\eta =T_{q}(\M)$$
holds isomorphically with universal constants, for any $0<\eta<1$,
where $\frac1q=1-\frac\eta2$. The latter fact is proved by using
the ``square argument'' of \cite{pis-intHpI} and  the Riesz type
factorization for triangular matrices (see \cite{xu-studia} for
more details).
 \end{rk}

We can now proceed to the second inequality of \eqref{int-ineq}.
Let $x\in L_{p_0}(d_{0})\cap L_{p_1}(d_{1})$ with
$\|x\|_{L_p(d_{\th})}<1$. Thanks to Remark \ref{triproj}, we have
$$\|T_+(x)d_\th\|_p< C\max(p,\,p').$$
As usual, let $\Delta=\{z\in\C :  0< Re(z)< 1\}$. We now appeal to
Lemma~\ref{pisierinter}. By  basic facts on complex interpolation
we find a continuous function $F: \overline \Delta \to
T_{p_0}(\M)\cap T_{p_1}(\M)$, which is holomorphic in $\Delta$ and
such that
 \begin{eqnarray*}
 F(\theta)= T_+(x)d_\th\quad\mbox{and}\quad
 \sup_{t\in \R}\big\{ \,||F({\rm i}t)\|_{p_0},\; ||F(1+{\rm i}t)\|_{p_1}\,\big\}
 \le Ct_p\max(p,\,p')\;{\mathop =^{\rm def}}\;C_p\,.
 \end{eqnarray*}
Put $G(z)=F(z)d_0^{z-1}d_1^{-z}$. Since $d$ is discrete and
bounded with bounded inverse, $G$ takes its values in
$T_{p_0}(\M)\cap T_{p_1}(\M)$, is continuous on $\overline \Delta$
and holomorphic in $\Delta$. We have $G(\theta)=T_+(x)$. On the
other hand, by Lemma~\ref{compa}, for $t\in \R$
 $$\| G({\rm i}t)\|_{L_{p_0}(d_0)}=
 \|T_+(G({\rm i}t))\|_{L_{p_0}(d_0)}\leq 2\|T_+(G({\rm i}t))d_0\|_{p_0}
 =2\|F({\rm i}t)d_0^{{\rm i}t}d_1^{-{\rm i}t}\|_{p_0}\le 2 C_p\,.$$
Similarly,
 $$\| G(1+{\rm i}t)\|_{L_{p_1}(d_1)}\le2 C_p\,.$$
It follows that
 $$\|T_+(x)\|_{(L_{p_0}(d_0),L_{p_0}(d_1))_\theta}\le2C_p\,.$$
Arguing in the same way for $T_-(x)$, we get the second inequality
of \eqref{int-ineq}.

\medskip

As for the duality of $L_p$-spaces, the other inequality is obtained
by duality. Applying the second inequality of \eqref{int-ineq} to
$p_1', p_0'$ and $1-\th$ instead of $p_0, p_1$ and $\th$ respectively,
we see that the identity
 $$\iota : L_{p'}(d_{1-\theta})\to(L_{p_1'}(d_0),\;
 L_{p_0'}(d_1))_{1-\theta}$$
is  bounded. We will dualize this inclusion.  The difficulty here lies on the
identifications.

First, we reformulate the previous result in
terms of non-weighted $L_p$-spaces.  As $d=\sum_{i=1}^k \lambda_i
e_i$, the map $\Sigma_d=L_{d}+ R_{d}$ is  a bijection on $L_q(d)$,
for any $1\leq q\leq \infty$. By definition $\Sigma_d$ is an
isometry from $L_q(d)$ onto $L_q(\M)$. With this in mind, we can
view the compatible couple $(L_{p_1'}(d_0),\; L_{p_0'}(d_1))$ as
$(L_{p_1'}(\M),\;L_{p_0'}(\M))^{\rm t}$ via a twisted
identification coming from the map
 $${\rm t}=\Sigma_{d_0}^{-1}\Sigma_{d_1} :
 L_{p_1'}(\M) \to L_{p_0'}(\M).$$
Then the maps
 \begin{eqnarray*}
  V_0=\Sigma_{d_{1-\theta}}^{-1} \Sigma_{d_0}  :
  L_{p'}(\M)\to  L_{p_1'}(\M)\\
  V_1=\Sigma_{d_{1-\theta}}^{-1} \Sigma_{d_1}  :
  L_{p'}(\M)\to  L_{p_0'}(\M)
 \end{eqnarray*}
are compatible with respect to $(L_{p_1'}(\M),L_{p_0'}(\M))^{\rm
t}$ (i.e., ${\rm t}\circ V_0=V_1$), so by interpolation they extend to a bounded map
 $$ V : L_{p'}(\M) \to (L_{p_1'}(\M),\;
 L_{p_0'}(\M))^{\rm t}_{1-\theta}\,.$$
By  duality we find
 $$ V^* : (L_{p_0}(\M),\;L_{p_1}(\M))^{{\rm t}^*,\theta}=
 (L_{p_1}(\M),\;L_{p_0}(\M))^{{\rm t}^*,1-\theta}\to
 L_{p}(\M).$$
Since
 $$(X_0,\; X_1)_\th\subset (X_0,\; X_1)^\th $$
isometrically for any compatible couple $(X_0,\; X_1)$ of Banach
spaces (see \cite{bergh}), we can restrict $V^*$ to
$(L_{p_0}(\M),\;L_{p_1}(\M))^{{\rm t}^*}_{\theta}$:
 $$ V^* : (L_{p_0}(\M),\;L_{p_1}(\M))^{{\rm t}^*}_\theta\to
 L_{p}(\M).$$
On the other hand, by duality, the compatibility
 of the couple $(L_{p_0}(\M),\;L_{p_1}(\M))^{{\rm t}^*}$  comes
from the map
 $$\Sigma_{d_0}^{-1}\Sigma_{d_1} :   L_{p_0}(\M) \to
 L_{p_1}(\M).$$
This is due to the fact that all $\Sigma_{d_i}$ are selfadjoint on
$L_2(\M)$. Thus ${\rm t}^*={\rm t}$ formally.  Also note that
$V^*$ is the extension of the compatible maps
 \begin{eqnarray*}
 V_1^*=\Sigma_{d_{1-\theta}}^{-1} \Sigma_{d_1}  :
 L_{p_0}(\M)\to L_{p}(\M),\\
 V_0^*=\Sigma_{d_{1-\theta}}^{-1} \Sigma_{d_0} :
 L_{p_1}(\M)\to  L_{p}(\M).
 \end{eqnarray*}

Now we return back to  the compatible couple
$(L_{p_0}(d_0),\;L_{p_1}(d_1))$. Then note that  the maps
 $$
 \Sigma_{d_0} :  L_{p_0}(d_0) \to  L_{p_0}(\M)\quad\mbox{and}\quad
 \Sigma_{d_1} :  L_{p_1}(d_1) \to L_{p_1}(\M)
 $$ are compatible isometries.  Composing them with $V^*$, we get a
bounded map from the interpolated space
$(L_{p_0}(d_0),\;L_{p_1}(d_1))_\theta$ to $L_{p}(\M)$, which extends
the following compatible maps
 \begin{eqnarray*}
 \Sigma_{d_0}\Sigma_{d_1}\Sigma_{d_{1-\theta}}^{-1} :
 L_{p_0}(d_0) \to  L_{p}(\M),\\
 \Sigma_{d_0}\Sigma_{d_1} \Sigma_{d_{1-\theta}}^{-1} :
 L_{p_1}(d_1) \to  L_{p}(\M).
 \end{eqnarray*}
Next, composing the last resulting map with the isometry
$\Sigma_{d_{\theta}}^{-1} :  L_{p}(\M)\to L_{p}(d_{\theta})$,
we deduce that the map
 $$\Sigma_{d_0}\Sigma_{d_1}
 \Sigma_{d_\theta}^{-1}\Sigma_{d_{1-\theta}}^{-1} :
 \left(L_{p_0}(d_0),\;L_{p_1}(d_1)\right)_\theta \to
 L_{p}(d_{\theta})$$ is bounded. Namely, \beq\label{schur-int}
 \|\Sigma_{d_0}\Sigma_{d_1}
 \Sigma_{d_\theta}^{-1}\Sigma_{d_{1-\theta}}^{-1}(x)
 \|_{L_{p}(d_{\theta})}\le C_{p'}
 \|x\|_{\left(L_{p_0}(d_0),\;L_{p_1}(d_1)\right)_\theta}\,.  \eeq
 Finally, to get the first inequality of \eqref{int-ineq}, we then
 just need to correct the left hand side above using
 Corollary~\ref{correc} for
 $\Sigma_{d_0}^{-1}\Sigma_{d_1}^{-1}\Sigma_{d_\theta}\Sigma_{d_{1-\theta}}$
 corresponds to a bounded Schur multiplier.  Thus we obtain the first
 inequality of \eqref{int-ineq}. This finishes the proof of step~1.

\begin{rk}
 Alternately, we can also first prove the first inequality of
\eqref{int-ineq} as in the appendix of \cite{pis-khin}, which is
essentially an argument dual to the previous one. Then we deduce
the second inequality by duality as above.
 \end{rk}

\medskip\n{\it Step~2:  Approximation.}
 Let
 $$\M_d=\bigcup_{n\ge1} q_n (\M\cap L_1(\M)) q_n\,,$$
where $q_n=\chi_{[n^{-1},\, n]}(d)$. It is easy to check that
$\M_d$ is a dense subspace of $L_p(f(d))$, for $1\le p\le\8$
(relative to the w*-topology for $p=\8$) and for any nondecreasing
$f$ on $\R_+$. Note that, for any $x\in \M_d$, $x$ belongs to
$q_n\M q_n$ for some $n$. As $q_n$ commutes with $d$, for such $x$
we have $\|x\|_{L_q(\M,f(d))}=\|x\|_{L_q(q_n\M q_n,f(q_nd))}$. On
the other hand, it is clear that $L_q(q_n\M q_n,f(q_nd))$ is a
contractively complemented subspace of $L_q(\M,f(d))$. Thus, it is
enough to prove the assertion for the reduced algebra $q_n\M q_n$,
with $q_nd$ instead of $d$. Therefore,  we can assume that both
$d$ and $d^{-1}$ are bounded operators on $\M$. In this case, $\M_d= L_1(\M)\cap\M$.

\medskip

Now let $(d_n)$ be a sequence of invertible positive operators
with discrete spectrum in the von Neumann subalgebra generated by
$d$ such that
 $$\|f_i(d_n)-f_i(d)\|_\8\leq \frac1n\,.$$
(For instance, each $d_n$ can be a positive linear combination of
mutually orthogonal spectral projections of $d$.) Then, for any
$1\leq q\leq \infty $, $i=0,\, 1$, and $x\in L_1(\M)\cap\M$, we
have
 $$\lim_n \|x\|_{L_q(f_i(d_n))}=\|x\|_{L_q(f_i(d_n))}.$$
This is clear as
 $$\big|\|x\|_{L_q(f_i(d_n))}-\|x\|_{L_q(f_i(d))}\big|
 \leq\|(f_i(d_n)-f_i(d))x+x(f_i(d_n)-f_i(d))\|_q
 \leq \frac 2 n\, \|x\|_q\,.$$
We go to the interpolation space. Note that $L_1(\M)\cap\M$ is
dense in $(L_{p_0}(f_0(d)),\; L_{p_1}(f_1(d)))_\theta$.  Let $x\in
L_1(\M)\cap\M$ such that
 $$\|x\|_{(L_{p_0}(f_0(d)),\; L_{p_1}(f_1(d)))_\theta}<1.$$
Then by \cite[Lemma~4.2.3]{bl} there exists a function
 $$\Psi(z)=\sum_{k} \psi_k(z)x_k$$
such that
 \begin{eqnarray*}
 \Psi(\theta)= x\quad\textrm{and}\quad
 \sup_{t\in \R} \{ \,\|\Psi({\rm i}t)\|_{L_{p_0}(f_0(d))},\;
 \|\Psi(1+{\rm i}t)\|_{L_{p_1}({f_1(d)})}\,\} < 1,
 \end{eqnarray*}
where $(x_k)$ is a finite sequence in $L_1(\M)\cap\M$ and
$(\psi_k)$ a finite sequence  of continuous functions on
$\overline \Delta$, holomorphic in $\Delta$ and vanishing at
infinity. Using the same function $\Psi$, but for the couple
$(L_{p_0}(f_0(d_n)),\; L_{p_1}(f_1(d_n)))$, we deduce
 $$\|x\|_{(L_{p_0}(f_0(d)),\; L_{p_1}(f_1(d)))_\theta}\geq
 \lim_n \|x\|_{(L_{p_0}(f_0(d_n),\;
 L_{p_1}(f_1(d_n))_\theta}\,.$$
To get the converse inequality, we again use duality.
As $L_1(\M)\cap \M$ is dense in all $L_q(f_i(d))$, all what we
need to show is that, for any $y\in L_1(\M)\cap \M$,
 $$\lim_n \|y\|_{L_q(f_i(d_n))^*}=\|y\|_{L_q(f_i(d))^*}\,.$$
However,
 $$\|y\|_{L_q(f_i(d_n))^*}=\|(L_{f_i(d)}+R_{f_i(d)})^{-1}
y\|_{L_{q'}}$$ and
$$(L_{f_i(d)}+R_{f_i(d)})^{-1} y=\int_0^\infty
e^{-f_i(d)t}ye^{-f_i(d)t} {\rm d} t.$$ Thus the desired result follows
from the dominated convergence theorem because $y\in L_1(\M)\cap
\M$. Therefore, the proof of Theorem \ref{wi11} is complete.

\begin{rk}\label{os}
Theorem \ref{wi11} can not be extended to arbitrary positive functions $f_0$ and
$f_1$. Otherwise, we would have that the map $
\Sigma_{d_0}\Sigma_{d_1}\Sigma_{d_\theta}^{-1}\Sigma_{d_{1-\theta}}^{-1}$
is bounded on $L_{p}(d_{\theta})$ (see \eqref{schur-int}). In terms of Schur multipliers, this
would mean that
 $$\left(\frac{(\lambda_i+\lambda_j)(\mu_i+\mu_j)}
 {(\lambda_i^{1-\th}\mu_i^{\th}+\lambda_j^{1-\th}\mu_j^{\th})
 (\lambda_i^{\th}\mu_i^{1-\th}+\lambda_j^{\th}\mu_j^{1-\th})}\right)_{i,j}$$
is a bounded Schur multiplier on $S_p$, for any positive sequences
$(\lambda_i)$ and $(\mu_i)$, which is false (take $\lambda_i=1/\mu_i=i$).
\end{rk}

\bigskip


\begin{thebibliography}{79}

\bibitem{bergh}
J.~Bergh.
\newblock On the relation between the two complex methods of interpolation.
\newblock {\em Indiana Univ. Math. J.}, 28:775--778, 1979.

\bibitem{bl}
J.~Bergh and J.~L{\"o}fstr{\"o}m.
\newblock {\em Interpolation spaces.}
\newblock Springer-Verlag, Berlin, 1976.

\bibitem{cha}
D.C.~Champeney.
\newblock {\em A handbook of {F}ourier theorems.}
\newblock {\em Cambridge University Press}, 1987.


\bibitem{fk}
Th. Fack and H.~Kosaki.
\newblock Generalized {$s$}-numbers of {$\tau$}-measurable operators.
\newblock {\em Pacific J. Math.}, 123:269--300, 1986.

\bibitem{haag-Lp}
U.~Haagerup.
\newblock {$L\sp{p}$}-spaces associated with an arbitrary von {N}eumann
  algebra.
\newblock In {\em Alg\`ebres d'op\'erateurs et leurs applications en physique
  math\'ematique (Proc. Colloq., Marseille, 1977)}, volume 274 of {\em Colloq.
  Internat. CNRS}, pages 175--184. CNRS, Paris, 1979.

\bibitem{hjx}
U.~Haagerup and J.~Junge and Q.~Xu.
\newblock A reduction method for noncommutative $L_p$-spaces
and applications.
\newblock {\em Trans. Amer. Math. Soc.}, to appear.

\bibitem{holm}
T.~Holmstedt.
\newblock Interpolation of quasi-normed spaces.
\newblock {\em Math. Scand.}, 26:177--199, 1970.

\bibitem{jupar-ros}
M.~Junge and J.~Parcet.
\newblock Rosenthal's theorem for subspaces of noncommutative {$L\sb p$}.
\newblock {\em Duke Math. J.}, 141:75--122, 2008.

\bibitem{kos-int}
H.~Kosaki.
\newblock Applications of the complex interpolation method to a von {N}eumann
  algebra: noncommutative {$L\sp{p}$}-spaces.
\newblock {\em J. Funct. Anal.}, 56:29--78, 1984.

\bibitem{nelson}
Ed. Nelson.
\newblock Notes on non-commutative integration.
\newblock {\em J. Funct. Anal.}, 15:103--116, 1974.

\bibitem{pis-intHpI}
G.~Pisier.
\newblock Interpolation between {$H\sp p$} spaces and noncommutative
  generalizations. {I}.
\newblock {\em Pacific J. Math.}, 155:341--368, 1992.

\bibitem{pis-sim}
G.~Pisier.
\newblock {\em Similarity problems and completely bounded maps}, volume 1618 of
  {\em Lecture Notes in Mathematics}.
\newblock Springer-Verlag, Berlin, expanded edition, 2001.

\bibitem{pis-khin}
G.~Pisier.
\newblock Remarks on the noncommutative Khintchine inequalities for $0<p<2$.
\newblock {\em J. Funct. Anal.}, 256:4128--4161, 2009.

\bibitem{px-survey}
G.~Pisier and Q.~Xu.
\newblock Non-commutative {$L\sp p$}-spaces.
\newblock In {\em Handbook of the geometry of Banach spaces, Vol.\ 2}, pages
  1459--1517. North-Holland, Amsterdam, 2003.

\bibitem{terp}
M.~Terp.
\newblock $L_p$ spaces associated with von neumann algebras.
\newblock Notes, Math. Institute, Copenhagen Univ., 1981.

\bibitem{xu-studia}
Q.~Xu.
\newblock Applications du th\'eor\`eme de factorisation pour des fonctions \`a
  valeurs op\'erateurs.
\newblock {\em Studia Math.}, 95:273--292, 1990.

\end{thebibliography}

\bigskip\footnotesize{
\n Laboratoire de Mathématiques, Université de Franche-Comté,
25030 Besan\c{c}on Cedex,  France\\
eric.ricard@univ-fcomte.fr;\hskip.3cm quanhua.xu@univ-fcomte.fr}

\end{document}